# Constructing Double Magma with Commutation Operations
Charles C. Edmunds

*Abstract.*

A *double magma* is a nonempty set with two binary operations satisfying the interchange law. We call a double magma *proper* if the two operations are distinct and *commutative* if the operations are commutative. A *double semigroup* is a double magma for which both operations are associative. Given a group G we define a double magma $(G,*,\bullet)$ with the commutator operations $x * y = [x, y](= x^{-1}y^{-1}xy)$ and $x \bullet y = [y, x]$. We show that $(G,*,\bullet)$ is a double magma if and only if G satisfies the commutator laws $[x, y; x, z] = 1$ and $[w, x; y, z]^2 = 1$. Note that the law $[x, y; x, z] = 1$ defines the variety of 3-metabelian groups. If both these laws hold in G, $(G,*,\bullet)$ is proper if and only if there exist $x, y \in G$ for which $[x, y]^2 \neq 1$. B.H. Neumann [6] has given an example of such a group which is not metabelian; thus the associated double magma is proper and produces an example with some complexity. The double magma $(G,*,\bullet)$ is a double semigroup if and only if $G$ is nilpotent of class 2. In this case, $(G,*,\bullet)$ is a proper double semigroup if and only if there exist $x, y \in G$ for which $[x, y]^2 \neq 1$. We construct a specific example letting G be the dihedral group of order 16. In addition we comment on a similar construction for rings using Lie commutators.

**Introduction.**

In [3] Kock introduces the notion of a double semigroup in relation to two-fold monoidal categories. He proves that both cancellative double semigroups and inverse double semigroups must be commutative. He begins with a discussion of the Ekmann-Hilton argument [2] for double moniods. They show that any unitary double magma is necessarily a commutative double semigroup with the two operations identical. In this note we will discuss a natural construction of double magma based on a group with operations defined in terms of commutation. We will determine exact conditions on the group so that our construction yields (i) a proper double magma (Corollary 1.7) and (ii) a proper double semigroup (Corollary 1.8). We begin with a few definitions.

**Definition:** A *double magma* is a triple $(M,*,\bullet)$ for which M is a nonempty set, $*$ and $\bullet$ are binary operations on M, and M satisfies the *interchange law*,
$$(w*x)\bullet(y*z) = (w\bullet y)*(x\bullet z),$$
for all $w, x, y, z \in M$.

We will call a double magma a *proper* double magma if its operations are distinct and *improper* if its operations are identical. Our goal here will be to construct classes of examples of proper double magma. We refer to a double magma as *commutative, associative,* or *unitary* when both of its operations are, respectively, commutative, associative, or unitary (*i.e.* have an identity element). A *double semigroup* is an associative double magma.

In [2] Ekman and Hilton, prove:

**Theorem:** If $(M,*,\bullet)$ is a unitary double magma, then $(M,*,\bullet)$ is an improper, commutative double semigroup with the same identity elements for $*$ and $\bullet$.

In light of this theorem, if we are to produce proper examples, we must be sure that not both magma have identity elements. It is possible for one of the two magma to have an identity element and still produce a proper example. Letting $D = \{a,b\}$ and defining the operations as:

| $*$ | $a$ | $b$ |
|---|---|---|
| $a$ | $a$ | $a$ |
| $b$ | $a$ | $b$ |

| $\bullet$ | $a$ | $b$ |
|---|---|---|
| $a$ | $a$ | $a$ |
| $b$ | $a$ | $a$ |

it is easy to check that $(D,*,\bullet)$ forms a proper, noncommutative double magma, in fact a double semigroup. Probably the most natural example of a double magma in elementary mathematics is the set of integers with addition and subtraction. Note that since subtraction is not associative, this is a proper (noncommutative) double magma but not a double semigroup. Similarly natural is the set of nonzero rationals with multiplication and division.

The fundamental idea of our constructions is to begin with any group G and a 2-variable word $W(a,b) \in F_2$, the free group of rank two freely generated by $a$ and $b$. We define two operations on G by
$$x*y = W(x,y) \text{ and } x \bullet y = W(y,x)$$

for all $x, y \in G$. To make $(G,*,\bullet)$ a double magma, we must impose the interchange law in this context. The law $(w*x)\bullet(y*z) = (w\bullet y)*(x\bullet z)$ translates as

$$W(W(y,z)),W(w,x)) = W(W(y,w),W(x,z)).$$

Depending on the word $W$, the translated version of the interchange law determines the properites G must have to produce a double magma with this construction. If $(G,*,\bullet)$ is improper then the law
$$W(x,y) = W(y,x)$$

holds for all $x, y \in G$. Thus, knowing that there is a pair of elements $x_0, y_0 \in G$ such that $W(x_0,y_0) \neq W(y_0,x_0)$, we are ensured that our double magma is proper.

As a simple example of how our construction works, suppose we take $W(a,b) = ab^{-1}$. Then the interchange law $(w*x)\bullet(y*z) = (w\bullet y)*(x\bullet z)$ becomes $yz^{-1}xw^{-1} = yw^{-1}xz^{-1}$. Letting $x = y = 1$ we obtain commutativity of G. Note that if G is commutative the interchange law holds; thus we conclude that $(G,*,\bullet)$ is a double magma if and only if G is an abelian group. To determine if this is a proper magma, we must investigate the law $W(x,y) = W(y,x)$. From $xy^{-1} = yx^{-1}$ we deduce that $(xy^{-1})^2 = 1$ and, letting $y = 1$, we see that $G$ must be of exponent 2. To summarize, we have shown that for $W(a,b) = ab^{-1}$, $(G,*,\bullet)$ is a double magma if and only if G is abelian and this is proper exactly when $G$ is not of exponent 2. Assuming that $G$ is an abelian group and not of exponent 2, we note that, although $G$ is abelian, the double magma is not. The commutativity of either of the two operations, would imply that $G$ is of exponent 2. Thus if we assume that $G$ is not of exponent 2, we know that neither $*$ nor $\bullet$ is commutative. Thus we have a class of examples, one for each abelian group not of exponent 2, yielding proper double magma. The result on associativity is that if either operation is associative, then $G$ is of exponent 2. Thus this construction will not produce a proper double semigroup. To give a very specific example, we can select $G$ to be the cyclic group of order three, $C_3 = \langle a; a^3 = 1 \rangle$. In this case the operation tables generated by the construction are as follows.

| * | 1 | $a$ | $a^2$ |
|---|---|---|---|
| 1 | 1 | $a^2$ | $a$ |
| $a$ | $a$ | 1 | $a^2$ |
| $a^2$ | $a^2$ | $a$ | 1 |

| $\bullet$ | 1 | $a$ | $a^2$ |
|---|---|---|---|
| 1 | 1 | $a$ | $a^2$ |
| $a$ | $a^2$ | 1 | $a$ |
| $a^2$ | $a$ | $a^2$ | 1 |

The remainder of this note will be an investigation of the case in which $W(a,b) = a^{-1}b^{-1}ab$, the commutator of $a$ and $b$. Here we will produce nontrivial examples of proper double magma and proper double semigroups.

**Preliminaries.**
Let $G$ be a group. For $x, y \in G$, we represent the action of $y$ on $x$ by conjugation as $x^y = y^{-1}xy$ and $x^{-y} = (x^{-1})^y = (x^y)^{-1}$. The *commutator of x and y* is defined as $[x,y] = x^{-1}y^{-1}xy$. We denote $[[x_1,x_2],x_3]$ by $[x_1,x_2,x_3]$ and, more generally, $[[x_1,x_2,\ldots,x_{n-1}],x_n]$ by $[x_1,x_2,\ldots,x_n]$. Further, we denote the commutator $[[w,x],[y,z]]$ by $[w,x;y,z]$. The normal closure of $\{[x_1,x_2,\ldots,x_n]: \text{for each } i \ (2 \leq i \leq n),\ x_i \in G\}$ is the $n^{th}$ *term of the lower central series*, denoted $\gamma_n(G)$ $(n \geq 2)$ with $\gamma_1(G) = G$. If $\gamma_{n+1} = \{1\}$, we say $G$ is *nilpotent of class n*. Note that $\gamma_2(G) = G'$, the *derived group* of $G$. The *second derived group* of $G$, denoted $G''$, is the normal closure of the set $\{[w,x;y,z]: w,x,y,z \in G\}$. If $G'' = \{1\}$ we refer to $G$ as *metabelian*. $G$ is *3-metabelian* if all of its three generator subgroups are metabelian. In [5] I. D. Macdonald showed that the single law $[x,y;x,z] = 1$ defines the variety (equational class) of 3-metabelian groups precisely. Bachmuth and Lewin [1] proved that the law $[x, y, z][y, z, x][z, x, y] = 1$ also defines the variety of 3-metabelian groups.

We will make use of the following identities as needed.
 (I i) $[x,y] = x^{-1}x^y = y^{-x}y$
 (I ii) $[x,y]^{-1} = [y,x]$

(I iii) $\left[a^{-1},b\right]=[a,b]^{-a^{-1}}$ and $\left[a,b^{-1}\right]^{\wedge}=[a,b]^{-b^{-1}}$

(I iv) $[ab,c]=[a,c]^{b}[b,c]=[a,c][a,c,b][b,c]$

(I v) $[a,bc]=[a,c][a,b]^{c}=[a,c][a,b][a,b,c]$

**Construction 1:** We define two binary operations on a group $G$. For each $x,y \in G$, let
$$x*y=[x,y] \text{ and } x \bullet y=[y,x].$$

**Proposition 1.1:** In any group $G$ the following statements are equivalent;
  (i) $(G,*)$ is commutative,
  (ii) $(G,\bullet)$ is commutative,
  (iii) for each $x,y \in G$, $x*y = x \bullet y$
  (iv) for each $x,y \in G$, $[x,y]^2 = 1$.

*Proof:* We will establish that each of the first three statements is equivalent to the fourth. Note first that if $*$ is commutative, then for each $x,y \in G$ we have $x*y = y*x$. This translates into the commutator equation $[x,y]=[y,x]$. By identity (I ii) we have
$$[x,y]=[y,x]=[x,y]^{-1}.$$
This is equivalent to the law $[x,y]^2 = 1$ in $G$. The equivalence of (ii) to (iv) and (iii) to (iv) are derived very similarly. ∎

**Proposition 1.2:** Each of $*$ and $\bullet$ is associative on $G$ if and only if $[x,y,z][y,z,x]=1$ for each $x,y,z \in G$.

*Proof:* Associativity of $*$, that is, $(x*y)*z = x*(y*z)$ translates into commutators as $[[x,y],z]=[x,[y,z]]$. By identity (i) we have $[[x,y],z]=[[y,z],x]^{-1}$. Thus $[x,y,z]=[y,z,x]^{-1}$ and it follows that $[x,y,z][y,z,x]=1$. Thus the law $[x,y,z][y,z,x]=1$ holds in G if and only if $*$ is associative. The same result follows in the same manner for $\bullet$. ∎

Note that in both these propositions the conditions being investigated for the operations turn out to be "varietal" in the group $G$. Something close to this will occur for the interchange law, $(w*x)\bullet(y*z) = (w \bullet y)*(x \bullet z)$, which translates into the commutator law $[y,z;w,x] = [y,w;z,x]$. Written with the variables permuted for aesthetic reasons, we will call the group law

(CI)   $$[w,x;y,z]=[w,y;x,z]$$

the *commutator interchange law*. If this holds for a system $(G,*,\bullet)$, the impact on the structure of $G$ is interesting. We will show first that (CI) implies that G is 3-metabelian. But (CI) is not equivalent to 3-metabelian. We will identify the class of groups determined by (CI) precisely in Theorem 1.6. This will be proved by commutator calculations, for which we will need to establish some preliminary lemmas.

**Lemma 1.3:** In any group $G$, the following laws are equivalent;

(3M i) $[x,y;x,z]=1$,

(3M ii) $[x,y;y,z]=1$,

(3M iii) $\left[x,y;[x,z]^u\right]=1$.

*Proof:* We will show first that (M3 i) and (M3 ii) are equivalent. Then we will show that (M3 i) and (M3 iii) are equivalent.

Firstly, by identities (I ii) and (I iii) we have $[x,y;x,z]=\left[[y,x]^{-1};x,z\right]=[y,x;x,z]^{-[y,x]^{-1}}$.

Inverting and conjugating both sides by $[y,x]$, we obtain $[x,y;x,z]^{-[y,x]}=[y,x;x,z]$. Thus if either law holds, it implies the other.

Secondly, note that $[x,y;x,z]=1$ follows by letting $u=1$ in $\left[x,y;[x,z]^u\right]=1$. We must then show that (3M ii) implies (3M iii). Suppose then that we have $[x,y;x,z]=1$. Substituting $zu$ for $z$, we get $1=[x,y;x,zu]$. Applying identity (I v) twice we obtain

$$1=[x,y;x,zu]=\left[x,y;[x,u][x,z]^u\right]=\left[x,y;[x,z]^u\right][x,y;x,u]^{[x,z]^u}.$$

Note that by (3M i) $[x,y;x,u]^{[x,z]^u}=1^{[x,z]^u}=1$. Thus we have derived (3M iii). ∎

The labels on these three laws remind the reader of I.D. Macdonald's result that each of these laws defines the 3-metabelian variety.

**Lemma 1.4:** If G is 3-metabelian, then the following laws hold in G;

$(L_1)$ $[x,y,z;x,u]=1$

$(L_2)$ $[x,y;x,u,v]=1$

$(L_3)$ $[x,y,z;x,u,v]=1$

*Proof:* To establish $(L_2)$ we begin with ( 3M i) in the form $[x,y;x,u]=1$ and substitute the product $yz$ for $y$. Applying (I iv) twice, we obtain,

$$1=[x,yz;x,u]=\left[[x,z][x,y][x,y,z];x,u\right]=[x,z;x,u]^{[x,y][x,y,z]}[x,y;x,u]^{[x,y,z]}[x,y,z;x,u].$$

By (3M i) the first two factors are trivial, therefore our result follows.

Given $(L_1)$, we can derive $(L_2)$ as an equivalent law as follows.

$$1=[x,y,z;x,u]=[x,u;x,y,z]^{-1}, \text{ therefore } [x,u;x,y,z]=1.$$

Lastly, we derive $(L_3)$ from $(L_1)$.

$$1=[x,y,z;x,uv]=[x,y,z;x,v][x,y,z;x,u][x,y,z;x,u,v].$$

Since the first two factors are instances of $(L_1)$, we obtain $[x,y,z;x,u,v]=1$. ∎

**Lemma 1.5:** If $G$ is 3-metabelian (*i.e.* satisfies the law $[x,y;x,z]=1$), then the law $[w,x;y,z][w,y;x,z]=1$ holds in G.

*Proof:* Assuming $G$ is 3-metabelian, we start from our variant of Macdonald's law, (3M *ii*) $[x,y;y,z]=1$, and substitute $yw$ for $y$. Thus by identities (I *iv*) and (I *v*),

$$1 = [x, yw; yw, z] = \big[[x,w][x,y][x,y,w]; [y,z][y,z,w][w,z]\big].$$

We organize the calculation by letting $X = [y,z][y,z,w][w,z]$. Therefore, by (I *iv*) and (I *v*), we have

$$1 = \big[[x,w][x,y][x,y,w]; X\big] = [x,w;X]^{[x,y][x,y,w]}[x,y;X]^{[x,y,w]}[x,y,w;X].$$

We will consider each of these three factors separately: call them A, B, and C.

*Consideration of* A.
$$A = [x,w;X]^{[x,y][x,y,w]} = \big[x,w;[y,z][y,z,w][w,z]\big]^{[x,y][x,y,w]}$$
$$= [x,w;w,z]^{[x,y][x,y,w]}\big[x,w;[y,z][y,z,w]\big]^{[w,z][x,y][x,y,w]}.$$

The first factor is trivial by (3M *ii*), therefore,
$$A = \big[x,w;[y,z][y,z,w]\big]^{[w,z][x,y][x,y,w]} = [x,w;y,z,w]^{[w,z][x,y][x,y,w]}[x,w;y,z]^{[y,z,w][w,z][x,y][x,y,w]}$$

$$= \Big[x,w;\big[w,[y,z]\big]^{-1}\Big]^{[w,z][x,y][x,y,w]}[x,w;y,z]^{[y,z,w][w,z][x,y][x,y,w]}$$

$$= \big[x,w;w,[y,z]\big]^{-[w,[y,z]]^{-1}[w,z][x,y][x,y,w]}[x,w;y,z]^{[y,z,w][w,z][x,y][x,y,w]}.$$

Again, the first factor is trivial by (3M *ii*), thus $A = [x,w;y,z]^{[y,z,w][w,z][x,y][x,y,w]}$.

We will now argue that each of these four conjugating elements commutes with $[x,w;y,z]$.

$$\big[[x,w;y,z],[y,z,w]\big] = \Big[[y,z;x,w]^{-1},[y,z,w]\Big] = \big[[y,z;x,w],[y,z,w]\big]^{-[y,z;x,w]^{-1}} = 1$$

applying $(L_3)$ to the last commutator. Next we consider,

$$\big[[x,w;y,z],[w,z]\big] = \Big[\big[[w,x]^{-1};y,z\big],[w,z]\Big] = \Big[[w,x;y,z]^{-[w,x]^{-1}},[w,z]\Big]$$

$$= \Big[[w,x;y,z]^{-1},[w,z]^{[w,x]}\Big]^{[w,x]^{-1}} = \Big[[w,x;y,z]^{-1},[w,z]^{[w,x]}\Big]^{[w,x]^{-1}}.$$

Since $[w,z;w,x]=1$, by (3M *i*), we conclude that $[w,z]^{[w,x]} = [w,z]$, and hence we have

$$\big[[x,w;y,z],[w,z]\big] = \Big[[w,x;y,z]^{-1},[w,z]\Big]^{[w,x]^{-1}} = \big[[w,x;y,z],[w,z]\big]^{-[w,x;y,z]^{-1}[w,x]^{-1}} = 1$$

by $(L_1)$. It also follows by $(L_1)$ that $\big[[x,w;y,z],[x,y]\big]=1$. And, finally, $\big[[x,w;y,z],[x,y,w]\big]=1$ by $(L_3)$. In summary, we have shown that $A = [x,w;y,z]$.

*Consideration of* B:

$$B = [x,y;X]^{[x,y,w]} = [x,y;[y,z][y,z,w][w,z]]^{[x,y,w]}$$
$$= [x,y;[w,z]]^{[x,y,w]}[x,y;[y,z][y,z,w]]^{[w,z][x,y,w]}$$
$$= [x,y;[w,z]]^{[x,y,w]}[x,y;y,z,w]^{[w,z][x,y,w]}[x,y;y,z]^{[y,z,w][w,z][x,y,w]}$$

The last two commutators are trivial by $(L_2)$ and (3M *ii*), respectively. Also note that $[[x,y;w,z],[x,y,w]] = 1$ by $(L_3)$, thus the conjugate on the first commutator in the last line commutes with the base commutator. Thus we have derived $B = [x,y;w,z]$.

*Consideration of* C:

$$C = [x,y,w;X] = [x,y,w;[y,z][y,z,w][w,z]] = [x,y,w;w,z][x,y,w;[y,z][y,z,w]]^{[w,z]}.$$

Note that by (3M *ii*) the first factor is trivial. Thus we have
$$C = [x,y,w;[y,z][y,z,w]]^{[w,z]} = [x,y,w;y,z,w]^{[w,z]}[x,y,w;y,z]^{[y,z,w][w,z]}.$$
Let us refer to these factors as D and E.
$$D = [x,y,w;y,z,w]^{[w,z]} = [[y,x]^{-1},w;y,z,w]^{[w,z]} = [y,x,w;y,z,w]^{-[y,x]^{-1}[w,z]}.$$
It follows that $D = 1$ by $(L_3)$.
$$E = [x,y,w;y,z]^{[y,z,w][w,z]} = [[y,x]^{-1},w;y,z]^{[y,z,w][w,z]} = [y,x,w;y,z]^{-[y,x]^{-1}[y,z,w][w,z]}.$$
It follows that $E = 1$ by $(L_1)$. Thus $C = 1$.

Returning to the original issue, we now have $1 = ABC = [x,w;y,z][x,y;w,z]$. ∎

We are now in a position to state and prove our main result pertaining to this construction.

**Theorem 1.6:** The commutator interchange law holds in a group if and only the group is 3-metabelian with every commutator of the form $[w,x;y,z]$ either trivial or of order 2. This means that in varieties of groups the law $[w,x;y,z] = [w,y;x,z]$ is logically equivalent to the union of the laws $[x,y;x,z] = 1$ and $[w,x;y,z]^2 = 1$.

*Proof:* Suppose first that $[w,x;y,z] = [w,y;x,z]$ holds in a group $G$. Replacing $w$ by $x$ we obtain $1 = [x,x;y,z] = [x,y;x,z]$, thus the interchange law implies that $G$ is 3-metabelian. Since $G$ is 3-metabelian, Lemma 1.5 indicates that $[w,x;y,z][w,y;x,z] = 1$ holds in $G$. Combining this with the interchange law we have
$$1 = [w,x;y,z][w,y;x,z] = [w,x;y,z]^2.$$

As a result we conclude that the interchange law implies both of the laws as stated in the theorem.

Conversely, Lemma 5.1 shows that if G is 3-metabelian, then $[w,x;y,z][w,y;x,z]=1$ in G. Thus we have $[w,x;y,z]=[w,y;x,z]^{-1}$. Since our hypothesis implies that $[w,y;x,z]$ is of order two, we conclude that $[w,x;y,z]=[w,y;x,z]^{-1}=[w,y;x,z]$. Thus the interchange law holds in G. ∎

We can now combine the information we have gathered to reflect on what kinds of examples of double objects we can construct from groups using the operations of left and right commutation.

**Corollary 1.7:** Given a group $G$, $(G,*,\bullet)$ is a proper double magma if and only if G is a nonabelian, 3-metabelian group with derived group not of exponent 2 and second derived group satisfying the identity $[w,x;y,z]^2=1$.
*Proof:* This follows immediately from Proposition 1.1. ∎

**Corollary 1.8:** Given a group *G*,
    (i) $(G,*,\bullet)$ is a double semigroup if and only if G is nilpotent of class 2 and
    (ii) $(G,*,\bullet)$ is a proper double semigroup if and only if *G* is a nonabelian class two group with derived group not of exponent two.

*Proof:* Suppose first that $(G,*,\bullet)$ is a double semigroup. According to Theorem 1.6, for $(G,*,\bullet)$ to be a double magma, G must be 3-metabelian and satisfy the law $[w,x;y,z]^2=1$. Recall Bachmuth and Lewin prove in [1], that the Jacobi identity, [*x, y, z*][*y, z, x*][*z, x, y*] = 1, defines the variety of 3-metabelian groups. Since our hypothesis implies that *G* is 3-metabelian, we know that the Jacobi identity holds in *G*. Proposition 1.2 states that for $(G,*,\bullet)$ to be associative, it is necessary and sufficient that $[x,y,z][y,z,x]=1$. Combining this with the Jacobi identity, it follows that $[z,x,y]=1$. This law defines the variety of groups which are nilpotent of class two. Conversely, if *G* is of class two, then all commutators of weights three and higher are trivial; therefore, G is 3-metabelian and has trivial second derived group.

Suppose now that $(G,*,\bullet)$ is a proper double semigroup. Being a double semigroup, we know that *G* is a class two group. Being proper we know that its derived group is not of exponent two, thus is it nonabelian. Conversely, if *G* is a non-abelian class two group, then G satisfies the interchange law (which is an equation between commutators of weight four, and hence trivial). G being of class two makes $(G,*,\bullet)$ associative, and this is a proper double semigroup exactly when there is a commutator whose square is nontrivial. Since G is of class two, the derived group is abelian and the commutators in *G* generate it. If each commutator were of order two or less, then the derived group would be abelian or of exponent 2. ∎

We should note that the double systems constructed here all have a zero element, the identity element of *G*. The most general example of a proper double magma constructed in this manner would be based on the relatively free group of the subvariety of 3-metabelian groups determined by the law $[w, x; y, z]^2 = 1$. In [6], B.H. Neumann gives an example of a 3-metabelian which is not metabelian. His group is of order $2^{14}$, has a derived group not of exponent 2, and satisfies the identity $[w, x; y, z]^2 = 1$. Thus, if a nonmetabelian example were desirable, Neumann's group could be used to construct a proper double magma while taking full advantage of the solvability length of 3.

Practically speaking, to construct a proper double magma in our manner, G could be chosen to be any metabelian group, or a group which is nilpotent of class 3, as long as the square of some commutator is nontrivial. This could be realized, among other ways, by letting G be a finite metacyclic group of odd order. Alternately, one could select any dihedral group of order not equal to 1, 2, 4, and 8. For example if we were to select the dihedral group of order six,
$$D_3 = \langle a, b; a^3 = 1, b^2 = 1, ba = a^2 b \rangle,$$
the group is metabelian and $[a, b]^2 = a^2 \neq 1$. Thus the double magma $(D, *, \bullet)$ is a proper, noncommutative magma.

To construct a double semigroup, we require that $G$ be nilpotent of class two with derived group not of exponent two. In the most general case, we could take G to be the relatively free group of class two. A simple example would be any nonabelian group of order $p^3$, with $p$ an odd prime. These groups are of class two and have no elements of even order. For a group containing 2-elements we could select $D_8 = \langle a, b; a^8 = 1, b^2 = 1, ba = a^7 b \rangle$, the dihedral group of order 16. This is a metabelian group and $[a, b]^2 = (a^6)^2 = a^4 \neq 1$. Thus $(D_8, *, \bullet)$ is a proper noncommutative double semigroup. To be completely explicit, we give the Cayley table for $(D_8, *)$ below. The table for $(D_8, \bullet)$ is the table below with each entry replaced by its inverse in $D_8$.

| $\star$ | 1 | a | $a^2$ | $a^3$ | $a^4$ | $a^5$ | $a^6$ | $a^7$ | b | ab | $a^2 b$ | $a^3 b$ | $a^4 b$ | $a^5 b$ | $a^6 b$ | $a^7 b$ |
|---|---|---|---|---|---|---|---|---|---|---|---|---|---|---|---|---|
| 1 | 1 | 1 | 1 | 1 | 1 | 1 | 1 | 1 | 1 | 1 | 1 | 1 | 1 | 1 | 1 | 1 |
| a | 1 | 1 | 1 | 1 | 1 | 1 | 1 | 1 | $a^6$ | $a^6$ | $a^6$ | $a^6$ | $a^6$ | $a^6$ | $a^6$ | $a^6$ |
| $a^2$ | 1 | 1 | 1 | 1 | 1 | 1 | 1 | 1 | $a^4$ | $a^4$ | $a^4$ | $a^4$ | $a^4$ | $a^4$ | $a^4$ | $a^4$ |
| $a^3$ | 1 | 1 | 1 | 1 | 1 | 1 | 1 | 1 | $a^2$ | $a^2$ | $a^2$ | $a^2$ | $a^2$ | $a^2$ | $a^2$ | $a^2$ |
| $a^4$ | 1 | 1 | 1 | 1 | 1 | 1 | 1 | 1 | 1 | 1 | 1 | 1 | 1 | 1 | 1 | 1 |
| $a^5$ | 1 | 1 | 1 | 1 | 1 | 1 | 1 | 1 | $a^6$ | $a^6$ | $a^6$ | $a^6$ | $a^6$ | $a^6$ | $a^6$ | $a^6$ |
| $a^6$ | 1 | 1 | 1 | 1 | 1 | 1 | 1 | 1 | $a^4$ | $a^4$ | $a^4$ | $a^4$ | $a^4$ | $a^4$ | $a^4$ | $a^4$ |
| $a^7$ | 1 | 1 | 1 | 1 | 1 | 1 | 1 | 1 | $a^2$ | $a^2$ | $a^2$ | $a^2$ | $a^2$ | $a^2$ | $a^2$ | $a^2$ |
| b | 1 | $a^2$ | $a^4$ | $a^6$ | 1 | $a^2$ | $a^4$ | $a^6$ | 1 | $a^6$ | $a^4$ | $a^2$ | 1 | $a^6$ | $a^4$ | $a^2$ |
| ab | 1 | $a^2$ | $a^4$ | $a^6$ | 1 | $a^2$ | $a^4$ | $a^6$ | $a^2$ | 1 | $a^6$ | $a^4$ | $a^2$ | 1 | $a^6$ | $a^4$ |
| $a^2 b$ | 1 | $a^2$ | $a^4$ | $a^6$ | 1 | $a^2$ | $a^4$ | $a^6$ | $a^4$ | $a^2$ | 1 | $a^6$ | $a^4$ | $a^2$ | 1 | $a^6$ |
| $a^3 b$ | 1 | $a^2$ | $a^4$ | $a^6$ | 1 | $a^2$ | $a^4$ | $a^6$ | $a^6$ | $a^4$ | $a^2$ | 1 | $a^6$ | $a^4$ | $a^2$ | 1 |
| $a^4 b$ | 1 | $a^2$ | $a^4$ | $a^6$ | 1 | $a^2$ | $a^4$ | $a^6$ | 1 | $a^6$ | $a^4$ | $a^2$ | 1 | $a^6$ | $a^4$ | $a^2$ |
| $a^5 b$ | 1 | $a^2$ | $a^4$ | $a^6$ | 1 | $a^2$ | $a^4$ | $a^6$ | $a^2$ | 1 | $a^6$ | $a^4$ | $a^2$ | 1 | $a^6$ | $a^4$ |
| $a^6 b$ | 1 | $a^2$ | $a^4$ | $a^6$ | 1 | $a^2$ | $a^4$ | $a^6$ | $a^4$ | $a^2$ | 1 | $a^6$ | $a^4$ | $a^2$ | 1 | $a^6$ |
| $a^7 b$ | 1 | $a^2$ | $a^4$ | $a^6$ | 1 | $a^2$ | $a^4$ | $a^6$ | $a^6$ | $a^4$ | $a^2$ | 1 | $a^6$ | $a^4$ | $a^2$ | 1 |

If $(R,+,\times)$ is a ring, then we define the *ring commutator* (or *Lie commutator*) for each $x, y \in R$ as $\langle x, y \rangle = xy - yx$. We can construct proper double systems from rings as we did from groups.

**Construction 2:** We define two binary operations on a ring $R$ as follows. For each $x, y \in R$, let $x * y = \langle x, y \rangle$ and $x \bullet y = \langle x, y \rangle$.

The *ring commutator interchange law* takes the form

(RCI) $$\langle w, x; y, z \rangle = \langle w, y; x, z \rangle,$$

and calculations, much simpler than those for groups, show that $(R, *, \bullet)$ is a proper magma precisely when (RCI) holds in R and there exist $x, y \in R$ such that $2\langle x, y \rangle \neq 0$. Analogous to the group case, it is true that (RCI) is equivalent to the laws $\langle x, y; x, z \rangle = 0$ and $2\langle w, x; y, z \rangle = 0$ in equational classes of rings.

To obtain a proper double semigroup the commutator identity $\langle x, y, z \rangle = 0$ must hold in $R$ as well as the existence of $x, y \in R$ such that $2\langle x, y \rangle \neq 0$. The law (RCI) is a polynomial identity holding on $R$, and this fact has a significant impact on the structure of $R$. We ask if there is a nice, structural characterization of those rings in which (RCI) holds?

Mathematics Department, Mount Saint Vincent University, Halifax, Nova Scotia, B3M 2J6, Canada
Email: cedmunds@eastlink.ca